\documentclass[11pt]{article}
\usepackage{amsmath}
\usepackage{amsthm}
\usepackage{mathtools}
\usepackage{amssymb}

\usepackage{enumitem}
\usepackage{hyperref}
\hypersetup{pdfcreator=XeLaTeX with hyperref}

\usepackage{geometry}
\usepackage{setspace}
\title{\textbf{{\normalsize NON-INVARIANCE OF THE NEVANLINNA CLASS UNDER THE ACTION OF THE CES\'ARO AND VOLTERRA OPERATORS}}} 
\author{\textsc{Apollon G. Paraskevas} \\ \\ Department of Mathematics \\ National and Kapodistrian University of Athens}
\date{}

\begin{document}
\maketitle 
\begin{abstract}
\noindent
Using a result of Hayman, we show that the Nevanlinna class of holomorphic functions on the unit disc is not invariant under the action of the Ces\'aro operator and more generally under the action of Volterra operators, $T_g$, provided that $\dfrac{1}{g'}$ belongs to the Nevanlinna class. \\ 
\end{abstract} 
\noindent
\textit{AMS classification numbers:} 30H15, 47B38  \medskip \\ 

A function $f$ belongs to the Nevanlinna class if $f$ is holomorphic in the open unit disc and the following condition is satisfied: 
\begin{equation*}
\sup\limits_{r\in(0, 1)}\bigg\lbrace \int_{0}^{2\pi}\log^{+}\Big\vert f\left(r\,e^{i\theta} \right) \Big\vert\,d\theta \bigg\rbrace<\infty,
\end{equation*}
as it is written in [1], where $\log^{+}$ is defined as 
\[
\log^{+}x=\begin{dcases}
\log x, \,\, x\geq 1 \\
0, \,\, 0\leq x<1 
\end{dcases},\quad
x\in\mathbb{R}^{+}.
\]
It is known from a theorem of F. and R. Nevanlinna (see [1]) that $f\in N$ if, and only if, $f(z)=\dfrac{a(z)}{b(z)}$ for $a, b$ holomorphic and bounded in the open unit disc $D(0, 1)$, and $b(z)\neq 0$ for every $z\in D(0, 1)$. Using this theorem, one can easily show that $N$ is an algebra over the field $\mathbb{C}$. \\ \\ \indent
Let $C$ denote the Ces\'aro operator. When $C$ acts on a function that is holomorphic on the open unit disc, it has an analytic form, i.e. if $f(z)=\sum\limits_{n=0}^{\infty}a_n z^n$, then $$C(f)(z)=\sum\limits_{n=0}^{\infty}\dfrac{1}{n+1}\sum\limits_{k=0}^{n}a_k z^n$$ and an integral representation (as in [3]): 
\[
C(f)(z)=\dfrac{1}{z}\int_{[0, z]}\dfrac{f(\zeta)}{1-\zeta}\,d\zeta.
\]
The Ces\'aro operator belongs in the family of Volterra operators $T_g$, where $$T_g(f)(z)=\dfrac{1}{z}\int_{[0, z]}f(\zeta)g'(\zeta)\,d\zeta.$$ In particular, if $g(z)=\log\left(\dfrac{1}{1-z}\right)$, then $g'(z)=\dfrac{1}{1-z}$, hence $$T_g(f)(z)=\dfrac{1}{z}\int_{[0, z]}\dfrac{f(\zeta)}{1-\zeta}\,d\zeta.$$ In [2], Hayman showed the existence of a function $h\in N$ for which $$R(h)(z)=\int_{[0, z]}h(\zeta)\,d\zeta$$ does not belong in $N$, where $H(z)=R(h)(z)$ is the primitive of $h$, vanishing at the origin. Using this result and the properties of the Nevanlinna class, we are going to show the following: \\ \\ 
\textbf{Proposition} \textit{There exists a $g\in N$ such that $C(g)\not\in N$.} \medskip \\ \indent
\textit{Proof.} Let $h\in N$ be the function in [2] such that $$H(z)=\int_{[0, z]}h(\zeta)\,d\zeta$$ does not belong in $N$. We proceed by showing the argument by contradiction and, thus, we suppose that for every $g\in N$ we have $C(g)\in N$. If we define $w$ by $w(z)=z$, then $w$ is bounded and holomorphic in $D(0, 1)$, hence we have that $w\in N$. We also have that $k(z)=1-z$ belongs to the Nevanlinna class, since $1-z$ is bounded and holomorphic. Now, from the fact that $N$ is an algebra and by the hypothesis, we have that $h(z)\cdot k(z) = h(z)(1-z)\in N$, which implies that $$C(h\cdot k)=\dfrac{1}{z}\int_{[0, z]}h(\zeta)\,d\zeta \in N.$$ It follows that $w(z)\cdot C(h\cdot k)(z)\in N$, from which we derive $$\int_{[0, z]}h(\zeta)\, d\zeta=H(z)\in N,$$
which is a contradiction. \qed \\ \\ \indent 
We will now give a sufficient condition on the symbol $g$ for the existence of such a function $f$ in the Nevanlinna class, so that $T_g(f)\not\in N$.\\ \\
\textbf{Theorem} \textit{If $\dfrac{1}{g'(z)}\in N$, then there exists a function $f\in N$ such that $T_g(f)\not\in N$.} \medskip \\  
\textit{Proof.} Let $h\in N$ with $R(h)\not\in N$. We set $f(z)=h(z)\cdot\dfrac{1}{g'(z)}\in N$ and, thus
\begin{equation*}
T_g(f)(z)=\dfrac{1}{z}\int_{[0, z]}h(\zeta)\dfrac{1}{g'(\zeta)}g'(\zeta)\,d\zeta = \dfrac{1}{z}\int_{[0, z]}h(\zeta)\,d\zeta.
\end{equation*} 
Since $w(z)=z\in N$, if we assume $$T_g(f)=\dfrac{1}{z}\int_{[0, z]}h(\zeta)\,d\zeta \in N,$$ then $$w(z)\cdot \dfrac{1}{z}\int_{[0, z]}h(\zeta)\,d\zeta\in N,$$ which yields $$R(h)(z)=\int_{[0, z]}h(\zeta)\,d\zeta\in N,$$
that is a contradiction. Hence, $T_g(f)\not\in N$. \qed \\ \\ \indent
In particular, the above Theorem implies that if $g'$ is holomorphic on the open unit disc $D(0, 1)$ and satisfies Re$g'(z)>0$ for all $z\in D(0, 1)$, then the Nevanlinna class is not invariant under the action of $T_g$. This holds because the assumption Re$g'(z)>0$ yields $\dfrac{1}{g'}\in N$ ([1]). \\ \\
\noindent
\textit{\textbf{Acknowledgements:}} The author would like to express his gratitude towards Professors Vassili Nestoridis and Aristomenis G. Siskakis.

\noindent
\textsc{Apollon G. Paraskevas} : Department of Mathematics, National and Kapodistrian University of Athens, Panepistemiopolis 157-84, Athens, Greece.\medskip \\
\textit{Email:} apollonparask@gmail.com


\begin{thebibliography}{99}
\bibitem{PDuren} Duren L. Peter, \textit{Theory of $H^p$ spaces}, Dover publications; Unabridged edition (2000). %1
%
\bibitem{Haym} Hayman W. K., \textit{On the characteristic of functions meromorphic in the unit disk and of their integrals}, Acta Math., Volume 112 (1964), 181-214.
%
\bibitem{Sisk} Siskakis G. Aristomenis, \textit{Composition Semigroups and the Cesaro Operator On $H^p$}, J. London Math. Soc. (2) 36 (1987), 153-164.
%
\end{thebibliography}
\end{document}